\documentclass[11pt,a4paper]{article}
\usepackage[scale=0.7,centering]{geometry}
\usepackage{amsfonts,amsthm}

\newtheorem{prop}{Proposition}
\newtheorem{thm}{Theorem}

\newtheorem{lemma}{Lemma}

\theoremstyle{definition}
\newtheorem{rmk}{Remark}

\newcommand{\enne}{\mathbb{N}}
\newcommand{\erre}{\mathbb{R}}
\newcommand{\Z}{\mathbb{Z}}

\newcommand{\E}{\mathbb{E}}
\renewcommand{\P}{\mathbb{P}}

\newcommand{\sgn}{\mathop{\mathrm{sgn}}\nolimits}
\renewcommand{\epsilon}{\varepsilon}

\newcommand{\ds}{\displaystyle}

\title{A class of stochastic games with infinitely many interacting
  agents related to Glauber dynamics on random graphs}
\author{Emilio De Santis\\
  \small\em Dipartimento di Matematica ``G. Castelnuovo''\\
  \small\em Universit\`a di Roma ``La Sapienza''\\
  \small\em P.le Aldo Moro 2, 00185 Roma, Italy\\
  \small\em e-mail: desantis@mat.uniroma1.it\\[12pt]
  Carlo Marinelli\\
  \small\em Institut f\"ur Angewandte Mathematik\\
  \small\em Rheinische Friedrich-Wilhelms-Universit\"at Bonn\\
  \small\em Wegelerstr. 6, D-53115 Bonn, Germany\\
  \small\em e-mail: cm788@uni-bonn.de\\[18pt]
}

\date{\normalsize July 6, 2007}

\begin{document}
\maketitle

\begin{abstract}\noindent
  We introduce and study a class of infinite-horizon non-zero-sum
  non-cooperative stochastic games with infinitely many interacting
  agents using ideas of statistical mechanics.
  First we show, in the general case of asymmetric interactions, the
  existence of a strategy that allows any player to eliminate losses
  after a finite random time. In the special case of symmetric
  interactions, we also prove that, as time goes to infinity, the game
  converges to a Nash equilibrium. Moreover, assuming that all agents
  adopt the same strategy, using arguments related to those leading to
  perfect simulation algorithms, spatial mixing and ergodicity are
  proved. In turn, ergodicity allows us to prove ``fixation'', i.e.
  that players will adopt a constant strategy after a
  finite time. The resulting dynamics is related to zero-temperature
  Glauber dynamics on random graphs of possibly infinite volume.
\end{abstract}

\section{Introduction}
The aim of this paper is to study a class of stochastic games with
infinitely many interacting agents that is closely connected with a
Glauber-type non-Markovian dynamics on random graphs. Let us briefly
explain the setting and our contributions both from the point of view
of game theory and of physics, referring to the next section for a
precise construction of the model. Our central results are theorems
\ref{t1}, \ref{teorema} and \ref{thm:fix} below.

We consider an infinite number of agents located on the vertices of
the two-dimensional lattice, where each agent is randomly linked with
others, and has positive or negative feelings regarding them.
Moreover, each agent is faced with the need of taking decisions that
affect himself and all others to whom he is linked. The objective of
each agent is to take (non-cooperative) decisions that ultimately do
not affect him negatively. Under a specific choice of the payoff
function of each player, we shall prove that there exists a decision
policy achieving this goal, and even more, that if each player adopts
this strategy a non-cooperative Nash equilibrium is reached.

From the physical point of view, we study a Glauber-type dynamics on a
random graph with the following features: the dynamics is
non-Markovian and has long-range interactions, in the sense that the
maximum distance between interacting particles is unbounded. For such
dynamics we prove spatial mixing (hence ergodicity) and fixation. To
the best of our knowledge, these problems are solved here for the
first time, even in the simpler case of a standard Glauber dynamics on
random graphs. Problems of dynamics on random graphs have attracted a
lot of attention in recent years (see the monograph \cite{Durrett} for
an extensive overview), as these structures are often more realistic
models of several phenomena than classical deterministic structures
(e.g. in network modeling, spread of epidemics, opinion formation,
etc.). For instance, C.~Cooper and A.~M.~ Frieze \cite{mini} prove the
existence of a critical coupling parameter at which the mixing time
for the Swendsen-Wang process on random graphs of sufficiently high
density is exponential in the number of vertices; M.~Dyer and
A.~Frieze \cite{mori} study the rapid mixing (in time) of Glauber
dynamics on random graphs with average degree satisfying a certain
condition (see also A.~Frieze and V.~Juan \cite{frosca} for a related
result).
In J.~P.~L.~Hatchett et al. \cite{hat}, the authors analyze the
dynamics of finitely connected disordered Ising spin models on random
connectivity graph, focusing on the thermodynamic limit.
I.~P\'erez Castillo and N.~S.~Skantzos \cite{hopp} study the Hopfield
model on a random graph in scaling regimes with finite average number
of connections per neuron and spin dynamics as in the Little-Hopfield
model.
On the other hand, as mentioned above, even though (spatial) mixing is
one of the most natural questions to ask about stochastic models of
interacting particle systems, it has not been discussed in the
literature, to the best of our knowledge. It is probably important to
recall that mixing is a key ingredient to obtain further results, such
as ergodicity. Moreover, just to cite another important application,
using Stein's methods (see e.g.  \cite{Billingsley}), mixing implies
the central limit theorem, which gives qualitative estimates on the
number of sites (or agents) with a positive spin (or opinion) in large
regions of the graph.

We would like to stress that our results on mixing are quite general,
and if one is only interested in the physical aspect of our work,
they could essentially skip the part of the paper which deals
with stochastic games, and concentrate only on the physical aspect.

Let us briefly discuss how the model and results of the present paper
are related to the existing literature on using methods of the theory
of interacting particle systems in economic modelling and game theory.
One of the first and still most cited works on the subject is a paper
by H.~F\"ollmer \cite{foellmer}, who considered a pure exchange
economy with (countably) infinitely many agents, each of which having
random preferences and endowments. In particular, agents are located
on the vertices of the $d$-dimensional lattice $\Z^d$, and their
preferences can be influenced by all his neighbors (i.e. such that
their euclidean distance is one). The author then considers the
problem of existence of a price system stabilizing the economy. See
also E.~Nummelin \cite{Nummelin} for further results in this
connection, but with a finite number of agents.
In U.~Horst and J.~Scheinkman \cite{HS} the authors study a system of
social interactions where agents are located on the nodes of a subset
of $\Z^d$, and each of them is provided with a utility function and a
set of feasible actions. The behavior of an agent is assumed to depend
on the choices of other agents in a reference group, which can be
random and unbounded. The authors, in analogy to our case, work under
the assumption that the probability of two agents being linked decays
with distance, and are concerned with the existence of equilibrium (in
the classical microeconomic sense).
U.~Horst \cite{horst-games} determines conditions such that
non-zero-sum discounted stochastic games with agents interacting
locally and weakly enough have a Nash equilibrium. While the set of
feasible actions in this paper is much richer than in ours, we do not
assume to have any knowledge on the reference group of each agent,
apart of being finite almost surely. We only allow agents to be able
to observe the dynamics of a (local) configuration around them. As a
result of the structural differences in the settings, the optimal
strategy in \cite{horst-games} is Markovian, while in our case it can
never be Markovian.

In general, the following features of our setting and results could be
particularly interesting from a game-theoretic perspective: we
consider games where interactions among agents are not known a priori,
and we explicitly construct a strategy that leads the game to
equilibrium, while the typical result of game theory is the existence
of equilibrium and a characterization of optimal strategies \emph{at}
equilibrium.

Let us also briefly recall that several other models of interacting
particle systems admit a natural interpretation in terms of social
interaction. Well-known examples are the voter process (see e.g.
T.~Liggett \cite{liggett}), used in models for the formation and
spread of opinions, or the Sherrington-Kirkpatrick model of spin
glasses (see e.g. section 2.1 in M.~Talagrand \cite{talagrand}).
Infinite interacting particle systems have found applications in
sociology as well: see, for instance, T.~Liggett and S.~Rolles
\cite{LR} for a model of formation of social networks.

The organization of the paper is as follows: in section \ref{sec:form}
we describe the model so we show how agents interact and what their
aim is; in section 3 a general strategy achieving the goal of each
agent is given, and section 4 proves spatial mixing, hence ergodicity,
of the dynamics, when all agents adopt the same strategy. Finally,
using the results on spatial mixing and ergodicity, we prove that the
game ``fixates'', i.e. that agents will adopt a constant strategy
after a finite time.  This phenomenon is reminiscent of the fixation
of zero-temperature dynamics (see e.g.  E.~De Santis and C.~M.~Newman
\cite{DSN}, L.~R.~Fontes \cite{fontes}, O.~H\"aggstr\"om \cite{hagg},
S.~Nanda, C.~M.~Newman and D.~L.~Stein \cite{NNS}).

\section{Model and problem formulation}\label{sec:form}
Let us first introduce some notation used throughout the paper. We consider
the two-dimensional lattice $\Z^2$ with sites $x=(x_1,x_2)$ and distance
$d$ defined by
\begin{equation}
\label{eq:distanza}
d(x,y) = |x_1-y_1| + |x_2-y_2|.
\end{equation}
The cardinality of a subset $\Gamma \subseteq \Z^2$ is denoted by
$|\Gamma|$.
We denote by $\Lambda_M$ the set of all $x\in\Z^2$ such that $d(O,x) \leq
M$, with $O=(0,0)$. If $x\in\Z^2$, $\Lambda_M(x)$ stands for $\Lambda_M+x$.
Our configuration space is $S=\{-1,+1\}^{\Z^2}$. The single spin space
$\{-1,+1\}$ is endowed with the discrete topology, and $S$ with the
corresponding product topology. Given $\eta \in S$, or equivalently
$\eta:\Z^2 \to \{-1,+1\}$, and $\Lambda \subseteq \Z^2$, we denote by
$\eta_\Lambda$ the restriction of $\eta$ to $\Lambda$.
Given a graph $G=(V,E)$, where $V$ and $E$ are the sets of its vertices and
edges, respectively, we shall denote by $\{x,y\}$ an element of $E$
connecting $x$, $y\in V$. 
For any $x\in V$, we shall denote by $\rho_x$ the distance of the longest
edge having $x$ as endpoint, namely we define
$$
\rho_x = \sup_{y:\,\{x,y\}\in E} d(x,y).
$$
Recall that the distance in variation of two probability measures
$\mu$ and $\nu$ on a discrete set $\Omega$ is defined as
$$
\|\mu-\nu\| = \frac12 \sum_{\omega\in\Omega} |\mu(\omega)-\nu(\omega)|.
$$

\smallskip

We shall now introduce an idealized model of a large ensemble of
interacting individuals. The ingredients will be a random graph, a
function on its edges (specifying an environment, roughly speaking),
and a stochastic process with values in $S$ describing the time
evolution of the system.\\
Let $\mathcal{G}=(V,E)$ be a random graph, whose set of vertices $V$
is given by all sites of the 2-dimensional lattice $\Z^2$, and whose
set of edges $E$ satisfy the following conditions: edges exist with
probability one between each site $x$ and all $y$ such that
$d(x,y)=1$, and
\begin{equation}
\label{eq:link-finiti}
\P(|\{y:\,\{x,y\}\in E\}|<\infty)=1 \quad \forall x\in V.
\end{equation}
We suppose that each site is occupied by an individual (we shall often
identify individuals with the sites they occupy, when no confusion will
arise), and that relations among individuals are modeled by the edges
of $\mathcal{G}$ and by a function $ j:V \times V\to\{-1,0,+1\}$; 
$j(x,y)= 0 $ if $\{ x,y \} \notin E$, otherwise $j(x,y) \in \{-1,+1\}  $. 

In particular, we shall say that individuals $x$ and $y$ are linked if
$\{x,y\}\in E$, and the value $j(x,y)$ shall account for the
``feelings'' of $x$ towards $y$: we set $j(x,y)=+1$ if $x$ is a
``friend'' of $y$, and $j(x,y)=-1$ if $x$ is an ``enemy'' of $y$.
We do not assume symmetry of $j$, i.e. friendship of an individual
towards another may not be reciprocal. Moreover, we assume that
individuals do not know with whom they are connected, nor whether
these individuals are friends or enemies.
Note also that in this model $x$ can be friend of $y$, $y$ friend of
$z$, but $x$ and $z$ can be either friends or enemies (a phenomenon
also called frustration in physics).

Let us now introduce a stochastic process $\sigma:[0,\infty) \to S$
modelling the evolution of the ``action'' (or opinion) of the
individuals. We shall use a graphical construction of the process, which
provides a specific version of basic coupling, i.e. it provides versions of
the whole family of stochastic processes on $\mathcal{G}$ (or on any finite
subset of it), all on the same probability space.
We assume that the initial configuration $\sigma_0$ is
chosen from a symmetric Bernoulli product measure. Moreover, the
continuous-time dynamics of $\sigma_t$ is given by independent Poisson
processes (with rate 1) at each site $x \in V$ corresponding to
those times $(t_{x,n})_{n\in\mathbb{N}}$ when the individual $x$ is asked
to update his opinion. Before describing the set of feasible ways of
opinion updating, let us introduce a reward for a generic individual $x$
at time $t_{x,n}$, as a result of his action:
$$
h_t(x) = \sgn
\left(\sum_{y:\,\{x,y\}\in E} j(x,y)\sigma_t(x)\sigma_t(y) \right),
$$
where we have set, for simplicity, $t\equiv t_{x,n}$.

We allow $x$ to base his decision on the history of
$\sigma_{\Xi_s}(s)$, $s\geq 0$, and $h(x)$, where $\Xi_s$ are
finite balls centered in $x$ with random radius which is nondecreasing
with respect to $s$, finite almost surely for all $s\geq 0$, and not 
'exploding'.
Formally, the decision of individual $x$ at time $t_{x,n}$ is a
$\{-1,+1\}$-valued random variable $u_{x,n}$ measurable with respect
to the $\sigma$-algebra generated by $\{\sigma_{\Xi_s}(s),\;s\leq
t_{x,n}\}$ and $\{h_s(x),\;s\leq t_{x,n}\}$, where $\Xi_s$ are
balls centered in $x$ such that
$$
\Xi_\infty = \lim_{s\to\infty} \Xi_s
$$
exists and is finite with probability one. We shall denote by
$\mathcal{E}^x_t$ the filtration just defined.

The dynamics of $\sigma$ is then completely specified by the updating rule
$$
\sigma_{t_{x,n}}(x) = u_{x,n}.
$$
Several remarks are in order: the reward $h_t(x)$ obtained by
individual $x$ as a result of his decision at time $t=t_{x,n}$ is
positive if the difference between pleased and damaged friends is
bigger than the difference between pleased and damaged enemies, negative
if the opposite happens, and zero if the value is the same. Since at
a fixed arrival time $t=t_{x,n}$ of the Poisson clock of $x$ no other
clock is ringing, i.e. $\P(t_{y,m}=t)=0$ for all $y\neq x$ and for all
$m\in\enne$, the dynamics of $\sigma$ is well-defined (also using the
graphical construction).  Finally, at any positive time $t$,
$\sigma_t(x)$ represents the last decision taken by individual $x$ up
to time $t$.

We formulate the following problem for the generic individual $x$:
find a strategy $\pi_x=(u_{x,1},u_{x,2},\ldots)$ such
that 
$$
h_t(x) \geq 0 \qquad \textrm{a.s.}
$$
for all $t\geq T_x$, where $T_x$ is a finite (random) time.

\begin{rmk}
  We built the random graph $\mathcal{G}$ on the two-dimensional
  lattice $\Z^2$ to give a ``geographic'' dimension to the problem and
  to have a simple notion of distance on the graph. However, all
  results in the next section still hold replacing $\Z^2$ with any
  higher dimensional lattice $\Z^d$, $d\geq 3$. We shall see below
  that choosing $d=2$ also affects a constant appearing in an
  assumption used to prove spatial mixing.
\end{rmk}

\section{Admissible strategies that eliminate losses}
In this section we construct explicitly a strategy $\pi_x$ for the
generic individual $x$ that asymptotically eliminate negative rewards,
i.e. such that $\P(h_t(x) \geq 0)=1$ for all $t$ greater than a random
time, which is finite with probability one. It will also be clear that
this strategy is non-cooperative, that is $\pi_x$ eventually eliminate
negative rewards irrespectively of the strategies adopted by all other
individuals.

For simplicity of notation let us describe the strategy
$\pi\equiv\pi_0$ for the individual located at the origin $O$. The
arrival times of his Poisson process and the corresponding decisions
and rewards will be denoted by $t_n$, $u_n$, and $h_n$, $n\in\enne$,
respectively.

The strategy $\pi=(u_1,u_2,\ldots)$ is best defined algorithmically
through a decision tree. We also need an additional ``data
structure'', i.e. a collection $\mathcal{R}$ of ordered triples of the
type $(\eta,u,h)$, where $\eta \in S$ is supported on finite balls,
$u\in\{-1,+1\}$, and $h\in\{-1,0,+1\}$.

At the first arrival time $t_1$, $u_1$ is chosen accordingly to a
Bernoulli law with parameter 1/2 (a ``fair coin toss''), and
$(\sigma_{\Lambda_1},u_1,h_1)$ is added to $\mathcal{R}$.

The description of the algorithm then follows inductively: at time
$t_{n+1}$, let $\Lambda_{M_n}$ be the support of the last configuration
added to $\mathcal{R}$.
Let $\sigma':=\sigma_{\Lambda_{M_n}}(t_{n+1}-)$ and check whether there exists
$(\sigma',u',h')\in\mathcal{R}$.
\begin{itemize}
\item If yes, set
$$
u_{n+1} = u' \frac{h'}{|h'|},
$$
with the convention $0/|0|:=1$. The reward $h_{n+1}$ corresponding to
$u_{n+1}$ is now obtained.
\begin{itemize}
\item If $h_{n+1}\geq 0$, no further action is needed.
\item If $h_{n+1}< 0$, then add to $\mathcal{R}$ the triplet
$(\sigma_{\Lambda_{M_n+1}},u_{n+1},h_{n+1})$.
\end{itemize}
\item Otherwise, set $u_{n+1}=u_n$, and add to
$\mathcal{R}$ the triplet $(\sigma',u_{n+1},h_{n+1})$.
\end{itemize}
The above algorithm formalizes the following heuristic procedure: the
agent starts looking at the configuration on the smallest ball
centered around him and plays tossing a coin. The next time his clock
rings, he checks whether he has already seen such a configuration. If
it is a new one, he will again memorize it and play by tossing a coin,
while if it is a known one he will play as he did before if he got a
positive reward, or the opposite way if he got a negative reward. Of
course it could happen that this way of playing still does not
guarantee a positive reward, in which case he will memorize the
configuration on a larger ball around himself and its associated
outcome.
\begin{rmk}
  One of the key steps of the algorithm requires one to look for a
  triplet $(\sigma',u',h')$ in $\mathcal{R}$, given
  $\sigma'=\sigma_{\Lambda}$, for a certain $\Lambda \subset \Z^2$.
  This operation is uniquely determined, i.e. there can exist only one
  triplet $(\sigma',u',h')\in\mathcal{R}$ with a given $\sigma'$. This
  can be seen as a consequence of the structure of the algorithm
  itself. Namely, as soon as the player ``observes'' the same
  configuration $\sigma'=\sigma_\Lambda$ with a different associated
  outcome $h$, he will immediately enlarge the support of observed
  configurations $\Lambda$.
\end{rmk}

We shall now prove that the strategy just defined eliminates losses
for large times.
\begin{thm} \label{t1}
For any individual $x$ there exists a random time $T_x$, finite with
probability one, such that
$$
\P\Big( h_{t\vee T_x}(x)\geq 0 \Big) = 1.
$$
\end{thm}
\begin{proof}
Let us define a sequence of random times $(\tau_n)_{n\in\enne}$ as follows:
$$
\tau_n = \inf \{ n\in\enne \,|\,
\exists (\overline{\sigma},u,h)\in\mathcal{R},\,
\mathrm{supp}\,\overline{\sigma}=\Lambda_{n+1}
\},
$$
with the convention that $\inf\emptyset=+\infty$.  In other words,
$\tau_n$ is the first time that individual $x$ includes into his
information set $\mathcal{R}$ the box $\Lambda_{n+1}$ (and if this never
happens, then $\tau_n=+\infty$).  Let $\tau_k$ be the last finite
element of the sequence $(\tau_n)_{n\in\enne}$. By assumption
(\ref{eq:link-finiti}) we know that $|\{y:\,\{x,y\}\in E\}|$ is finite, hence 
$k \leq \rho_x$ because $h_t(x)$ only depends on those $y$ linked to $x$, for
all times $t$. Therefore the biggest $\Lambda_n $ observed by the agent 
in the origin is finite.
\par\noindent 
Define the family of sets
$$
A_k(t) = \left\{
\sigma_{\Lambda_k}(t_{x,\ell}):\, t_{x,\ell} \in [\tau_k,t]
\right\}.
$$
It clearly holds $A_k(t_1) \subseteq A_k(t_2)$ for $t_1<t_2$,
hence we can define
$$
A_k(\infty) = \lim_{t\uparrow\infty} A_k(t).
$$
Since $A_k(t) \subset \{-1,+1\}^{\Lambda_k}$, and $\Lambda_k$ is finite, then
there exists $T_x>0$ such that $A_k(T_x) = A_k(\infty)$, hence
$$
A_k(t) = A_k(\infty) \qquad \forall t>T_x.
$$
We claim that $h_{x,n}\geq 0$ for all $t_{x,n}>T_x$. In fact, for
every $t_{x,n}>T_x$ there exists
$(\overline{\sigma},u,h)\in\mathcal{R}$ with
$\overline{\sigma}(y)=\sigma_{t_{x,n}}(y)$ for all $y\in\Lambda_k$.
But since $\tau_{k+1} = \infty$, the algorithm will give as output a
$u_n$ such that $h_n\geq 0$ (to convince oneself it is enough to
``run'' the algorithm). In a more suggestive way, one could say that
after $T_x$ individual $x$ has already been faced at least once with
all possible configurations that are relevant for him, and therefore
knows how to take the right decision.
\end{proof}
\begin{rmk}
  (i) Note that the strategies of other individuals never enter into the
  arguments used in the proof. Therefore individual $x$ is sure to
  reach the goal of eliminating losses in finite time irrespectively
  of the strategies played by all other individuals. 

  (ii) However, we would like to stress that the random time $T_x$ is
  not a stopping time (i.e. it is not adapted to the filtration
  $\mathcal{E}^x_t$). In fact, $T_x$ depends in general on the
  decisions of other individuals, whose policies are not necessarily
  adapted to $\mathcal{E}^x_t$. In general, even if all policies were
  adapted, the random times $\{T_x\}_{x \in \mathbb{Z}^2}$ would not
  be stopping times.

  (iii) Let us also observe that although we formally allowed the
  strategy $\pi_x$ to be adapted to $\mathcal{E}^x_t$, the information
  used by the strategy constructed in the proof of Theorem \ref{t1} is
  much smaller.  Similarly, one could refine the way the memory
  structure $\mathcal{R}$ is constructed, for instance by eliminating
  configurations on smaller balls, when one starts to add new
  configurations on balls of higher radius. However, we preferred to
  keep the construction of $\mathcal{R}$ as it is to avoid
  non-essential complications.
\end{rmk}
As a consequence of theorem \ref{t1} and of observation (i) in the
above remark, one has the following result, which essentially states
that the games admits an ``asymptotic'' Nash equilibrium.
\begin{prop}
  Let $M\in\mathbb{N}$ and assume that each player $x\in
  \Lambda_M:=[-M,M]\times[-M,M]$ adopts the strategy $\pi_x$ defined
  above.  Then there exists a finite random time $T_M$ after which no
  agent can gain by any change in their strategy given the strategies
  currently pursued by other players.
\end{prop}
It is important to observe that in the above proposition we implicitly
assume that each player only cares about ``not loosing'', or
equivalently he distinguishes only between ``loosing'' ($h_t(x)<0$)
and ``not loosing'' ($h_t(x)\geq 0$). In this sense, after $T_M$,
there is no point for any player $x\in\Lambda_M$ to change his
strategy, as proved in theorem \ref{t1}. The statement of the
proposition is in general false if the player distinguishes between
$h_t(x)>0$, $h_t(x)=0$, and $h_t(x)<0$.

We think that one can prove (and we leave it as a conjecture), that
this asymptotic equilibrium is not Pareto. This could be done adapting
ideas of O.~H\"aggstr\"om \cite{hagg}, who proved that
zero-temperature dynamics on a random graph does not reach the minimum
energy configuration.

\section{Spatial mixing and ergodicity}
The main result of this section, which plays an essential role in the
results about fixation of the next section, is that a spatial mixing
property holds. We shall work under the following hypothesis, which
states that the probability of two agents being linked decays
algebraically with their distance. 
\smallskip\par\noindent
\textbf{Standing assumption.\/} It holds that
\begin{equation}
\label{eq:edges}
\P(\{x,y\} \in E) \leq \frac{C}{d(x,y)^{9}},
\end{equation}
for all $y$ such that $d(x,y)>1$, where $C$ is a positive constant.
\smallskip\par\noindent
Note that assumption (\ref{eq:edges}) implies (\ref{eq:link-finiti}).
Moreover, the exponent appearing on the right-hand side of
(\ref{eq:edges}) depends on the dimension of the lattice and it
is needed in order to use well-known combinatorial estimates on path
counting in $\Z^2$ in the proofs to follow. However, it would not be
difficult to generalize our arguments to any higher dimensional
lattice $\Z^d$, $d\geq 3$, at the expense of replacing the exponent
$9$ with a (higher) constant depending on the dimension $d$, and of
using more complicated estimates in the proofs. Since this point is
not essential and would only add technical complications, we preferred
to fix $d=2$.

Before stating the main theorem of this section, we need to introduce
the following set of conditions.
\smallskip\par\noindent
\textbf{Hypothesis H.}
\textit{%
  The random graph $\mathcal{G}=(V,E)$ and the process
  $\sigma:[0,\infty) \to \{ -1,1 \}^{\Z^2}$,
  satisfy the following conditions:
\begin{enumerate}
\renewcommand{\labelenumi}{(\roman{enumi})}
\item For each vertex $x\in\Z^2$ there exists a Poisson process $P_x$,
  and the Poisson processes $\{P_x\}_{x \in \Z^2}$ are
  mutually independent. Denoting by $\Upsilon_x = \{t_{x,n}\}$ the set of
  arrival times of $P_x$, the value of $\sigma_t(x)$ is allowed to change
  only at times $t \in \Upsilon_x$.
\item Given any couple $(x,y)\in \Z^2 \times \Z^2$, the probability
  $\P(\{x,y\}\in E)$ is defined and it can depend on $d(x,y)$.
  Moreover, for any choice of $(x,y)$, $(v,w) \in \Z^2\times\Z^2$ with
  $(x,y) \neq (v,w)$, the events $\{x,y\}\in E$ and $\{v,w\}\in E$ are
  mutually independent.
\item The evolution of the process is local, i.e.
  $\sigma_{t_{x,n}}(x)$ is measurable with respect to
  $\mathfrak{F}^x_{t_{x,n}}$, where $\mathfrak{F}^x_t$ denotes the
  $\sigma$-algebra generated by $\{\sigma_s(y): \{x,y\} \in E \textrm{
    or } y=x, \; s<t\}$. We denote by $\mathcal{F}^V_t$ the
  $\sigma$-algebra generated by $\cup_{x\in V} \mathfrak{F}^x_t$.
\item Both the probability of two agents being linked and the
  evolution of $\sigma$ are translation invariant, i.e.  $\P(\{x,y\}
  \in E) = \P(\{x+v,y+v\} \in E)$ and $\P(\sigma_t\in
  A|\sigma_0=\eta)=\P(\sigma_t \in A+v|\sigma_0 (\cdot )=\eta(\cdot+v))$.
\end{enumerate}
}
\smallskip\par\noindent
We can now state the main theorem of this section.
\begin{thm}
\label{teorema}
If Hypothesis H holds true, then $\sigma$ satisfies the spatial mixing
property
\begin{equation}
\label{eq:teorema}
\lim_{\Lambda\to\Z^2}
\P\Big(\sigma_{\Lambda_0}(t)=\eta | \mathcal{F}^{\Lambda^c}_t\Big) =
\P(\sigma_{\Lambda_0}(t)=\eta),
\end{equation}
where $\Lambda_0$ is any finite region in $\Z^2$.
\end{thm}
Note that the process $\sigma$ is translation invariant if each agent
adopts the same strategy at each decision time (the strategy does
\emph{not} need to be the one defined in section 3).
Before giving the proof of the theorem, we establish some auxiliary
results.

We shall use the following terminology: by ``box of side length $L$''
we mean the set $[-L/2,L/2]^2 \subset \Z^2$. For $\rho<1$, we call
``subbox of side length $L^\rho$'' any one of the $L^{1-\rho}$ square
sets into which a box of side length $L$ can be subdivided. We always
assume $L^\rho$, $L^{1-\rho} \in \mathbb{N}$ (without loss of
generality, as it will be clear). Furthermore, we shall say that two
subboxes $R$ and $S$ are ``neighbors'' if $d(R,S)\leq\sqrt{2}$, so
every subbox has $8$ neighbor subboxes. We shall call ``path of
subboxes'' a sequence of subboxes $(R_k)_{k=1,\ldots,K}$ such that
$R_k$ and $R_{k+1}$ are neighbors for each $k=1,\ldots,K-1$. Two
subboxes $R$, $S$ are ``linked'' if there exist $x\in R$, $y\in S$
such that $\{x,y\}\in E$.

In the following lemma we introduce a sequence of boxes
increasing to $\Z^2$, each of one further subdivided into a variable
number of boxes also increasing to $\Z^2$, but at a lower rate.  

\begin{lemma}
\label{lem:scatole}
There exist a sequence of integer numbers $L_n \uparrow +\infty$, a
sequence of square boxes $Q_{L_n}$ of side length $L_n$, each of them
partitioned into subboxes of side $L_n^\rho$, $\rho=13/42$, such that
only a finite number of the boxes $Q_{L_n}$ will contain linked
non-neighbor subboxes.
\end{lemma}
\begin{proof}
We use a Borel-Cantelli argument on a suitable sequence of box side lengths
$L_n$.  In particular, let $L$ be a positive integer, $Q_L$ a square of
side $L$, subdivided into subboxes of side $L^\rho$.  The probability
of an agent $x$ to be linked with some other agent of a non-neighbor subbox
is bounded by
$$ 
\sum_{y:d(x,y)\geq L^\rho} \frac{C}{d(x,y)^9}
\leq
C_1 \int_{L^\rho}^{\infty} \frac{1}{v^9}\,2\pi v\,dv
=
C_2 \frac{1}{L^{7\rho}},
$$
where $C$, $C_1$, $C_2$ are positive constants.
Therefore two agents in non-neighbor subboxes exist with probability
not larger than
$$
L^2 \sum_{y:d(x,y)\geq L^\rho} \frac{C}{d(x,y)^9}
\leq
C_2 \frac{1}{L^{7\rho-2}} \to 0, \quad \textrm{as $L \to \infty$.}
$$
Taking now a subsequence $L_n$ growing to infinity rapidly enough,
$$
\sum \P(A_{L_n}) < \infty,
$$
where $A_{L_n}$ denotes the event that $Q_{L_n}$ contains linked
non-neighbor subboxes. By Borel-Cantelli lemma, only a finite number
of occurrences of $A_{L_n}$ can happen, which finishes the proof.
\end{proof}
Recall that for a sequence of i.i.d. standard exponential random
variables $\{X_i\}$ one has
\begin{equation}
\label{eq:rate}
\P\Big( \sum_{i=1}^n X_i < n\alpha \Big)
\leq
e^{-\Phi(\alpha)n}
\quad \forall \alpha < \E X_1,
\end{equation}
where the so-called rate function $\Phi$ is given by
$$
\Phi(\alpha) = \alpha-1-\log\alpha.
$$

\noindent
\emph{Proof of Theorem \ref{teorema}.}
We use a coupling argument to show that
\begin{equation}
\label{fultim}
\sup_{\zeta',\zeta''} \left|
\P(\sigma_{\Lambda_0}(t)=\eta | \sigma_{\Lambda^c}(0)=\zeta')
-
\P(\sigma_{\Lambda_0}(t)=\eta | \sigma_{\Lambda^c}(0)=\zeta'')
\right| \to 0
\end{equation}
and hence, by the inequality
$$
\begin{array}{l}
\displaystyle \sup_{\zeta',\zeta''} \left|
\P(\sigma_{\Lambda_0}(t)=\eta | \sigma_{\Lambda^c}(0)=\zeta')
-
\P(\sigma_{\Lambda_0}(t)=\eta | \sigma_{\Lambda^c}(0)=\zeta'')
\right|\\[10pt]
\quad \geq
\left|
\P(\sigma_{\Lambda_0}(t)=\eta | \sigma_{\Lambda^c}(0)=\zeta)
-
\P(\sigma_{\Lambda_0}(t)=\eta)
\right| \quad \forall \zeta,
\end{array}
$$
that (\ref{eq:teorema}) holds.\\
We construct two coupled systems $\sigma'$, $\sigma''$ on the same
probability space supporting $\sigma$ in the following way:
$\sigma_x'(0)=\sigma_x''(0)=\sigma_x(0)$ for all $x\in\Lambda$;
$\sigma'$ and $\sigma''$ update their state according to the same
translation-invariant rule of $\sigma$; all other randomness in the
system (the random graph, the Poisson processes, the ``coin tosses''
needed for the decision rules) coincide.
Define, for any $x\in V$, the random time
$$
\tau_x = \inf\{t\geq 0:\; \sigma'_x(t)\neq\sigma''_x(t)\},
$$
and introduce the process 
$$
[0,\infty) \times V \ni (t,x) \mapsto 
\nu_x(t) = 1(t\geq\tau_x) \in \{0,1\}.
$$
Using a pictorial language, we shall say that we color $x$ with
black as soon as the two processes $\sigma'$ and $\sigma''$ differ at
$x$.
Let us also introduce another process
$\tilde\nu:[0,\infty)\times V \to\{0,1\}$ with the property
$\tilde\nu_x(t) \geq \nu_x(t)$ a.s. for all $x$ and all $t$. The
dynamics of $\tilde\nu$ is specified as follows: $\tilde\nu_x(0)=0$
for all $x\in V$, and $\tilde\nu_x$ can turn to one as a consequence
of two classes of events. In particular, (i) $\tilde\nu_x(t)=1$ if
there exists $x'$ belonging to the same subbox of $x$ such that
$\nu_{x'}(t)=1$, and (ii) $\tilde\nu_x(\tau)=1$ if there exists $y$
belonging to a neighbor subbox such that $\tilde\nu_y(\tau)=1$, where
$\tau$ is any arrival time of the Poisson process relative to $x$.
Moreover, we assume that 1 is an absorbing state for $\tilde\nu_x$,
for all $x$. Using again a pictorial analogy, we could say that the
black area generated by $\tilde\nu$ is bigger than the black area
generated by $\nu$. In particular, as soon as a site $x$ turns black,
(i) implies that the whole subbox to which it belongs becomes black as
well.
\par\noindent
By Lemma \ref{lem:scatole}, there exists a positive integer $N$ and a
sequence $L_n$ such that for all $n>N$ the boxes $Q_{L_n}$ contain no
linked non-neighbor subboxes.
The shortest path of subboxes from the boundary of the box $Q_{L_n}$
to its center has length $L_n^{1-\rho}/2$ (therefore, for $n$ large
enough, the shortest path of subboxes from the boundary of the box
$Q_{L_n}$ to $\Lambda_0$ has length greater or equal than
$L_n^{1-\rho}/2-1$).
Setting $T_{Q_{L^\rho}}=\inf_{x\in Q_{L^\rho}} t_{x,1}$ (recall that
$t_{x,1}$ is the first arrival time of the Poisson processes relative
to $x$), one has that the distribution of $T_{Q_{L^\rho}}$ is
$\mathrm{Exp}(L^{2\rho})$, where $\mathrm{Exp}(\lambda)$ stands for
the law of an exponential random variable with parameter $\lambda$.
The minimum time for the formation of a path of $k$ ``black'' subboxes
along a fixed path (of sites) from the boundary of $Q_{L_n}$ to the
origin is given by
$$
T = \sum_{i=1}^k T_i
$$
for all $n>N$ (from now we shall tacitly assume $n>N$), where
$T_1,\ldots,T_k$ are i.i.d. exponential random variables with
parameter $L_n^{2\rho}$ (independence and the value $L_n^{2\rho}$
follow by the memoryless property of the exponential distribution).

Note that the sequence of subboxes in a path turning black does not
influence the minimum time needed for the formation of such path,
which is a sum of independent exponential random variables of
parameter $L_n^{2\rho}$, using again the memoryless property of
exponential distributions. It follows by (\ref{eq:rate}) that, for
$0<\alpha<1$, one has
$$
\P(T \leq k \alpha L_n^{-2\rho})
\leq
e^{-(\alpha-1-\log\alpha)k}.
$$
Denoting by $T_{\partial Q_{L_n} \to O}$ the (random) time needed to
form a path of black subboxes from the boundary of $Q_{L_n}$ to the
origin $O$, we obtain the estimate
\begin{eqnarray*}
\P\Big(
T_{\partial Q_{L_n} \to O} \leq
\frac{\alpha}{2}L_n^{1-3\rho}
\Big) 
&\leq&
4 L_n^{1-\rho} \sum_{k\geq \frac{L_n^{1-\rho}}{2}}
8 \cdot 7^{k-1} \exp\Big(-(\alpha-1-\log\alpha)k\Big),
\end{eqnarray*}
hence, for $0<\alpha<\frac{1}{7e}$,
$$
\lim_{n\to\infty} \P\Big(
T_{\partial Q_{L_n} \to O} \leq
\frac{\alpha}{2}L_n^{1-3\rho}
\Big) = 0.
$$
Here the term $4L_n^{1-\rho}$ accounts for the possible initial
subbox on the boundary of $Q_L$, and $8 \cdot 7^{k-1}$ is an upper
bound for the number of paths (of subboxes) of length $k$ starting in
a given subbox.
We obtain that, as $n\to\infty$, the term on the right hand side goes to
zero like $e^{-\beta L_n^{1-\rho}}$ (modulo polynomial terms), with
$\beta$ a positive constant.
Again by a Borel-Cantelli argument we obtain 
$$
\P\Big(
\lim_{L\to\infty} T_{\partial Q_L \to O} = \infty
\Big) = 1.
$$
Moreover, the evolution of the central subbox is completely independent
on the configuration outside $\Lambda$ until it turns black, and so
the theorem is proved.
\hfill$\Box$
%
\begin{rmk}
  Although (\ref{fultim}) has been proved only for a particular choice
  of a sequence of increasing boxes $\Lambda_n$, one can easily show
  that any increasing sequence of boxes will do. In fact, the supremum
  appearing in (\ref{fultim}) is decreasing with respect to $\Lambda$,
  hence it is enough to prove the theorem for any (fixed) subsequence.
\end{rmk}

\section{Fixation}\label{last}
In this section we shall work under the general assumptions introduced
in section 2 and 4, and furthermore we assume that each player
adopts the same strategy (hence the dynamics is translation
invariant), and that interactions are symmetric, i.e. that
$j(x,y)=j(y,x)$ for any $x$, $y \in V$. The latter hypothesis is
essential, as it would be possible to find counterexamples to our
results in the case of asymmetric interactions.
As before, we shall denote by $x$ an arbitrary agent, fixed throughout
this section. Let us define the random time $T_x$ as
\begin{equation}
\label{def:T_x}
T_x =
\sup \{t: \hbox{ at time $t$ agent $x$ sees a new configuration or loses} \}.
\end{equation}
As it follows from Theorem \ref{t1}, $T_x$ is finite with probability one.
Moreover, by definition, agent $x$ will not loose at any time after $T_x$.
Let us also define the random variable $M_x$ as the number of times
agent $x$ changes his state (i.e. updates his opinion) during the
time interval $(0,+\infty)$.

The main result of this section is the following:
\begin{thm}
\label{thm:fix}
Assume that each agent adopts the strategy constructed in section 3.
Then each agent $x\in V$ updates his opinion only a finite number of
times, i.e.
$$
\P(M_x<\infty) = 1.
$$
\end{thm}
Before proving theorem \ref{thm:fix}, we shall need some more
definitions and preparatory results.

Let us recall the definition of $\rho_x$:
\begin{equation}
\label{rr}
\rho_x = \sup_{y:\,\{x,y\}\in E} d(x,y),
\end{equation}
the distance from $x$ of his farthest connected agent.
Note that one has, as follows by the standing assumption (\ref{eq:edges}),
\begin{equation}
\label{mav}
\P(\rho_x \geq r) \leq \sum_{s \geq r } 4s \frac{C}{s^{8 +\epsilon
}}\leq \frac{K}{r^{6 +\epsilon }}
\end{equation}
where $C$ and $K$ are constants depending on $x$. Therefore
$\E\rho_x^k$, $1\leq k\leq 5$ are finite:
\begin{equation}
\label{ef}
\E\rho_x^k \leq 1 + \sum_{r=2}^{\infty} r^k \P(\rho=r) \leq 
1 + \sum_{r=2}^{\infty} r^k \P(\rho \geq r) \leq 
1 + \sum_{r=2}^{\infty} r^k \frac{K}{r^{6 +\epsilon }} < \infty.
\end{equation}

Let us also define the energy (or Lyapunov) function on a finite set
$\Lambda\subset\Z^2$ as
\begin{equation}
\label{Ham1}
H_\Lambda(\sigma) = - \sum_{u \in \Lambda} \tilde{h}_u(\sigma),
\end{equation}
where
\begin{equation}
\label{ham2}
\tilde{h}_u(\sigma) = \sum_{v: \{u,v\} \in E} j(u,v)\sigma_u\sigma_v.
\end{equation}
In the following we shall denote by $\Lambda_n$ the square box
$[-n ,n] \times [-n ,n]$.

\begin{lemma} \label{llee}
There exists a continuous function $e:\erre_+ \to [-\E\rho^2, \E\rho^2]$
such that
\begin{equation}
\label{gh}
\lim_{n\to \infty } \frac{H_{\Lambda_n}(\sigma(t))}{|\Lambda_n|}
= e(t) \quad \mathrm{a.s.}
\end{equation}
\end{lemma}
\begin{proof}
  By the definitions of
  $\tilde{h}_x(\sigma(t))$, $\rho_x$, it follows that for each time $t $
  $$
  -\rho_x^2 \leq \tilde{h}_x(\sigma(t)) \leq \rho_x^2,
  $$
  hence, taking expectations, recalling (\ref{ef}), and using
  translation invariance
  $$
 - \infty < -\E\rho_O^2 \leq \E\tilde{h}_x(\sigma(t)) \leq E\rho_O^2< \infty .
  $$
At any time $t$, using the space ergodicity of the system (implied
by the spatial mixing property proved in Theorem \ref{teorema}), we obtain
\begin{equation}
\label{nb}
\lim_{n\to\infty } \frac{H_{\Lambda_n} (\sigma(t))}{|\Lambda_n|} =
\E\tilde{h}_O(\sigma(t)) \quad \mathrm{a.s.}
\end{equation}
Setting $e(t)=\E\tilde{h}_O(\sigma(t))$, we just have to prove that $e$
is continuous. Using again the spatial ergodicity of $\sigma$, the
proportion of agents in $ \Lambda_n$ 
taking at least a decision in the time interval
$]t_1, t_2[$ tends to $1 - e^{-(t_2-t_1)} \leq t_2 - t_1$ as
$n\to\infty$.
Since each agent is endowed with a Poisson process that is independent
from all other processes and random variables describing the dynamics
of the system, the mean energy variation of each agent is bounded by
$\E\rho^2$. Therefore we also have
\begin{equation}
\begin{array}{rcl}
|e(t_2)-e(t_1)| &=& \ds \lim_{n\to\infty}
\frac{|H_{\Lambda_n}(\sigma(t_2))-H_{\Lambda_n}(\sigma(t_1))|}{|\Lambda_n|}\\[10pt]
&\leq&  (1- e^{-(t_2 -t_1 )}) \E\rho^2 \leq (t_2 -t_1 ) \E\rho^2,
\end{array}
\end{equation}
i.e. the function $e$ is Lipschitz continuous.
\end{proof}

Let us now define the following discrete random sets for agent $x$,
which are subsets of the set of arrival times of his Poisson process:
\begin{eqnarray*}
 N_1(x)  &=&  \{t : t \leq T_x,
\hbox{the agent in $x$ sees a known configuration at time $t$ and loses}\} \\
N_2 (x)   &=&      \{t : t \leq T_x,
\hbox{there is an arrival of the Poisson process in $x$}\}   \setminus   N_1(x)  \\
N_3(x)  &=&  \{t : t > T_x, \hbox{the agent in $x$ changes
opinion}\}
\end{eqnarray*}
Note that by definition of $T_x$, at any time $t>T_x$ agent $x$ can
only see known configurations, and can only win.
\par\noindent
We also define, for every $t>0 $ and $x \in \Z^2$, the random sets
$$
N_i(t,x) =  N_i(x) \cap [0,t],
$$
for $i=1, \,\,2,\,\, 3 $. 
\par\noindent
Moreover, for $\Lambda\subset\Z^2$, we set
$$
N_i(t,\Lambda) = \bigcup_{x\in\Lambda} N_i(t,x).
$$
\par\noindent
The dynamics of the system and the definition of $e(t)$ imply that
$e(t)$ is determined only by the changes of $\sigma_\tau(\cdot)$, $\tau \in
\{N_i(t,x)\}_{x\in\Z^2}$, $i=1, 2, 3$. We can therefore write
$$
e(t) = e_1(t) + e_2(t) + e_3(t) ,
$$
where $e_i(t)$ denotes the component of $e(t)$ determined by
changes of $\sigma_\tau(\cdot)$ for $\tau\in\{N_i(t,x)\}_{x\in\Z^2}$.
Moreover one has $e_2(t)\leq 0$ because we are eliminating the
arrivals where the agent lost, and in this case the energy can only
decrease.

We are now in the position to prove the theorem on the fixation of the
stochastic dynamics.

\smallskip\par\noindent
\emph{Proof of Theorem \ref{thm:fix}.}\/
In virtue of the translation invariance of the system, it is enough to
prove the result for the agent in the origin. First observe that $M_O
\leq |N_1(O)| +|N_2(O)| + |N_3(O)|$, because $N_1(O)$ and $N_2(O)$ may
contain Poisson arrival times in which the agent $O$ does not change
his opinion. Denoting $N_i(O)$ by $N_i$ for simplicity of notation, we
shall prove that $|N_i|<\infty$ almost surely for $i =1,2, 3$.
Let us first observe that the following inclusion relations hold:
$$
 \{ | N_1 (O)|+| N_2(O) | =\infty \} \subset  \{ T_O = \infty \} \cup (\bigcup_{n \geq 1}
  \{ N_1 ( O ) =\infty , \,\, T_O <n \}),
$$
and
$$
\begin{array}{l}
\ds \{|N_1(O)|+|N_2(O)| = \infty, \; T_O<n \} \subset \\
\ds \qquad\qquad
\{ |\textrm{Arrivals in $(0,n)$ of the Poisson process in the origin} |= \infty, \; T_O <n \}.
\end{array}
$$
Recalling that $\P(\{T_O = \infty\})=0$ we obtain
$$
\begin{array}{l}
\ds P(|N_1(O)|+|N_2(O)| = \infty ) \leq \\
\ds \qquad\qquad \sum_{n=1}^{\infty}\P(
|\textrm{Arrivals in $(0,n)$ of the Poisson process in the origin} | = \infty
) = 0.
\end{array}
$$
Thus we only need to show that $| N_3 (O) |$ is almost surely
finite.
First we observe that one has
\begin{eqnarray} \label{bn}
e_1(t) &=& \lim_{n\to\infty} \frac{1}{|\Lambda_n|}
\sum_{\tau\in N_1(t,\Lambda_n)}
H_{\Lambda_n}(\sigma(\tau))-H_{\Lambda_n}(\sigma(\tau^-)) \leq 
\E\rho_O^3 \quad \textrm{a.s.},
\end{eqnarray}
because the number of changes in the origin $N_1 (t,O)$ is at most
$\rho_O$ (the maximum number of enlargements of the box observed by
the agent $x$), and in any change the energy can increase at most by
$\rho_O^2$. Finally, the spatial ergodicity yields the almost sure
upper bound in (\ref{bn}).
\par\noindent
At any time $\tau\in N_3$ the energy $H_\Lambda(\sigma(t))$
decreases at least of one unit, i.e. $H_\Lambda(\sigma(t)) \leq
H_\Lambda(\sigma(t^-)) -1$, otherwise the agent does not change
opinion. Thus
\begin{equation} \label{ql}
\begin{array}{rcl}
\ds e_3(t) &=& \ds \lim_{n\to\infty} \frac{1}{|\Lambda_n|}
  \sum_{\tau\in N_3(t,\Lambda_n)}
     H_{\Lambda_n}(\sigma(\tau))-H_{\Lambda_n}(\sigma(\tau^-)) \\
\ds &\leq& \ds \lim_{n\to\infty}
\frac{-|N_3(t,\Lambda_n)|}{|\Lambda_n|} = -\E |N_3(t,O)|,
\end{array}
\end{equation}
where we have used once more the spatial ergodicity.

By Lemma \ref{llee} and noting that the energy is initially zero
(because agents choose $+1$ or $-1$ with probability $1/2$), one has the following inequality
$$
-\E\rho^2 \leq e(t) = e_1(t)+e_2 (t)+e_3(t)  \leq e_1 (t)+e_3(t),
$$
which holds uniformly in time $t$. Using inequalities (\ref{bn})
and (\ref{ql}) we obtain $\E|N_3(t,O)| \leq \E\rho_O^3 + \E\rho_O^2
\leq \infty$ uniformly in $t$, hence also in the limit as $t \to
\infty$. But $\E|N_3(O)|<\infty$ obviously implies $|N_3 (O)|<\infty$
a.s., so we have shown that $M_O \leq | N_1(O) | + | N_2(O) | + |
N_3(O) | < \infty $ and the proof is complete.\hfill $\Box$

\begin{rmk}
We can also deduce, following the proof of Theorem \ref{thm:fix},
that
$$
P( N_3 >C   )\leq \frac{ E ( \rho_O^2 ) +
 E ( \rho_O^3 )  }{C},
$$
as an immediate consequence of Markov's inequality.
\end{rmk}

\begin{rmk}
  Let us briefly comment on the connection between the fixation result
  just proved and the results of De Santis and Newman \cite{DSN}. The
  improvement is twofold: namely, the dynamics considered here does
  not coincide (locally) with zero-temperature dynamics. It is
  immediate to prove that at any given time there is at least an agent
  which does not follow the zero-temperature dynamics. This implies
  that on any time interval the zero-temperature dynamics and our
  dynamics are almost surely different. Our could say, perhaps
  somewhat informally, that our dynamics is a perturbation of
  zero-temperature dynamics with the property of preserving fixation.
  Moreover, as already mentioned several times, our dynamics is
  non-Markovian, while the arguments used in \cite{DSN} hold only for
  Markovian dynamics.
\end{rmk}

\end{document}